\documentclass[12pt]{article}
\usepackage{tikz-cd}
\usepackage{mathrsfs}
\usepackage{amsfonts}
\usepackage{amsthm,amsmath,amssymb,anysize}
\newtheorem{lemma}{Lemma}[section]
\newtheorem{theorem}[lemma]{Theorem}

\newtheorem{proposition}[lemma]{Proposition}
\newtheorem{definition}[lemma]{Definition}
\newtheorem{corollary}[lemma]{Corollary}
\newtheorem{example}[lemma]{Example}
\usepackage[numbers,sort&compress]{natbib}
\usepackage{amsmath}
\allowdisplaybreaks[3]
\marginsize{3cm}{3cm}{3.6cm}{3cm}
\numberwithin{equation}{section}

\title{\bf{The Universal Enveloping Algebras of Lie-Rinehart Color Algebras}}
\author{{Xiaohui Chi and Liming Tang\footnote{corresponding author}}\\
\small{School of Mathematical Sciences}\\
\small{Harbin Normal University}\\
\small{150025 Harbin, China}\\
\small{E-mail: limingtang@hrbnu.edu.cn}}
\date{ }
\date{ }
\begin{document}
\maketitle
 \begin{quotation}
{\small\noindent \textbf{Abstract}: In the present paper, Lie-Rinehart color algebras over a ring $\mathbb{R}$ are introduced, which combine Lie-Rinehart (super)algebras and Lie color algebras. Then we establish some classical algebraic concepts such as derivations, homomorphisms, and representations of Lie-Rinehart color algebras and obtain some examples. Furthermore, we describe the universal enveloping algebras for Lie-Rinehart color algebras, which play  key roles in the studies of representations and homology. Finally, the Poincar\'e-Birkhoff-Witt theorem is obtained under certain conditions of the Lie-Rinehart color algebras. 

 \vspace{0.05cm} \noindent{\textbf{Keywords}}: universal envelope algebra; Lie-Rinehart color algebra; $\Gamma$-graded associative algebra

\vspace{0.05cm} \noindent \textbf{Mathematics Subject Classification
2010}: 17B35, 17B75}
\end{quotation}
 \setcounter{section}{0}
\tableofcontents
\section{Introduction}
In this paper, our aim is to introduce Lie-Rinehart color algebras over a ring $\mathbb{R}$, 
with emphasis on basic results for further research on structures and representation.  
Besides being useful as a new algebraic system, Lie-Rinehart color algebras integrate the 
core structures of both Lie-Rinehart (super)algebras and  Lie color algebras, which play important 
roles in many parts of mathematics and physics. Lie-Rinehart algebras are the algebraic 
counterparts of Lie algebroids and provide a framework for studying
geometric structures \cite{ref33}, which also play key roles in Poisson geometry and in the 
study of Gerstenhaber and Batalin-Vilkovisky algebras
\cite{ref34,ref35}. Lie color algebras, originating from the study of supersymmetry and 
parastatistics in particle physics, are a natural generalization of Lie superalgebras 
\cite{ref36}.  


Lie-Rinehart algebras originate from the intersection of differential geometry and commutative algebra. 
 On a smooth manifold $M$, the smooth functions form a commutative associative algebra $C^\infty(M)$, while the smooth vector fields form a Lie algebra $\mathfrak X(M)$ under the Lie bracket. Vector fields act on smooth functions as derivations, and multiplication by smooth functions makes $\mathfrak X(M)$ a natural $C^\infty(M)$-module.
The early development of this notion can be traced back to Herz's Lie pseudo algebras \cite{ref1} and Palais's
$d$-Lie rings \cite{ref2}. Rinehart subsequently developed
the theory systematically and constructed the universal enveloping algebra \(U(A,\mathfrak g)\) of a Lie-Rinehart algebra \((A,\mathfrak g)\) and proved a Poincar\'e--Birkhoff--Witt theorem \cite{ref3}.
 
 Let us recall the definition of a Lie-Rinehart algebra \cite{ref3,ref15}.
 Let $A$ be a commutative $\mathbb{R}$ algebra. A Lie-Rinehart algebra $(A,\mathfrak{g},\rho)$ consists of a Lie algebra $\mathfrak{g}$ over $\mathbb{R}$, an $A$-module structure on $\mathfrak{g}$, and an anchor map $\rho:\mathfrak{g}\to \mathrm{Der}_{\mathbb{R}}(A)$ which is both an $A$-module homomorphism and a Lie algebra homomorphism, satisfying the relation:
\[
[x,ay]=a[x,y]+\rho(x)(a)y,\quad \mbox{for all}\,x,y\in \mathfrak{g},\ a\in A.
\]
 The Lie-Rinehart structure has also been generalized to the superalgebra setting. Roger \cite{ref19} generalized the theory of Lie-Rinehart algebras to the superalgebra framework and derived related conclusions for several classical examples of supermanifolds. Ehret and Makhlouf \cite{ref40} investigated low-dimensional classification and developed a formal deformation theory of Lie-Rinehart superalgebras. Lamkin \cite{ref4} studied the universal enveloping algebra of Lie-Rinehart superalgebras and established the corresponding PBW theorem. 
 Bouarroudj, Ehret, Makhlouf and Shyntas \cite{ref20} introduced restricted Lie-Rinehart superalgebras and constructed their universal enveloping algebras. Scholars have also successively studied various generalized forms of Lie Rinehart algebras. For example, 3-Lie-Rinehart superalgebras \cite{ref16}, graded 3-Lie-Rinehart algebras \cite{ref28}, graded Lie-Rinehart algebras \cite{ref14} and so on. Lie color algebras generalize Lie algebras and Lie superalgebras by introducing abelian group gradings and commutation factors. Early work on Lie color algebras includes the paper of Rittenberg and Wyler \cite{ref21}, and their structure was subsequently studied systematically by Scheunert \cite{ref22}.
Lie color algebras arise as a special case of the Lie-Rinehart color algebras introduced in this paper.

Motivated by \cite{ref3,ref19,ref21,ref32}, it is natural to introduce
$\Gamma$-graded Lie-Rinehart algebras with bicharacters, namely
Lie-Rinehart color algebras, and study their algebraic structures. We also describe the universal enveloping algebras and obtain some basic properties. Specifically, we replace the Lie algebra structure with the Lie color algebra structure, simultaneously endow the color-commutative algebra \(A\) and the module \( \mathfrak{g}\) with \(\Gamma\)-graded structures, and maintain the graded compatibility of the module action and the anchor map. 

The classification and deformation results in the super setting also suggest concrete questions for Lie-Rinehart color algebras. For fixed graded algebras \(A\) and Lie color algebras \(\mathfrak g\), one can ask which module actions and anchors give compatible Lie-Rinehart color structures, and how the grading and commutation factor affect their deformation theory. And the study of Lie-Rinehart superalgebras in characteristic \(p\) provides another direction for future research. 

The paper is organized as follows: In Section 2, we mainly review the basic concepts. In Section 3, we introduce the basic definition of Lie-Rinehart color algebras, provide corresponding examples, and introduce related concepts such as homomorphisms, modules, and ideals of Lie-Rinehart color algebras. In Section 4, we construct the universal enveloping algebra and study its basic properties. Drawing on Lamkin's work \cite{ref4},
we then prove a Poincar\'e--Birkhoff--Witt theorem under certain conditions. 

Throughout this paper, let $\mathbb{R}$ be an integral domain in which
$2$ and $3$ are invertible and \(\Gamma\) be an abelian group. All algebras and modules are over a ring $\mathbb{R}$,
and all gradings are by \(\Gamma\). let $\mathbb N$ denote the set of nonnegative integers. We write $\mathbb{R}^\times$ for the group of units
of $\mathbb{R}$. For a $\Gamma$-graded color-commutative algebra
$A=\bigoplus_{\alpha\in\Gamma} A_\alpha$, we simply write $A$
when its grading is understood.
\section{Preliminaries}
We recall some basic notions. For more details see \cite{ref27,ref30}. 
 Let \(\Gamma\) be an abelian group. A \textit{bicharacter} on \(\Gamma\) is a mapping \((\cdot, \cdot) : \Gamma \times \Gamma \to \mathbb{R}^\times\) such that
\[
(\alpha, \beta)(\beta, \alpha) = 1, \quad (\alpha, \beta + \gamma) = (\alpha, \beta)(\alpha, \gamma),
\]
\[
(\alpha + \beta, \gamma) = (\alpha, \gamma)(\beta, \gamma),
\]
for all \(\alpha, \beta, \gamma \in \Gamma\).
To avoid any confusion, we denote the zero element in \(\Gamma\) by \(\bar{0}\).
A \textit{color algebra} is just a \(\Gamma\)-graded associative \(\mathbb{R}\)-algebra. 
 Let $A,B$ be color algebras. A morphism \(f: A \to B\) of \(\mathbb{R}\)-algebras is a \textit{morphism of color algebras} if \(f\) is homogeneous of degree $1_{\Gamma}$ . A \textit{color ideal} is a \(\Gamma\)-graded two-sided ideal.
Let \(M\) be a \(\Gamma\)-graded $\mathbb{R}$-module. Then
\[
\text{{\rm{End}}}_\mathbb{R}(M)_\gamma := \{ f \in \text{{\rm{End}}}_\mathbb{R}(M) \mid f(M_\mu) \subseteq M_{\mu+\gamma} \}
\]
is the \(\mathbb{R}\)-submodule of homogeneous \(\mathbb{R}\)-linear maps of degree \(\gamma \in \Gamma\) of \(\text{{\rm{End}}}_\mathbb{R}(M)\) and
\[
\text{{\rm{End}}}_\mathbb{R}(M)_\Gamma := \bigoplus_{\gamma \in \Gamma} \text{{\rm{End}}}_\mathbb{R}(M)_\gamma
\]
is a color algebra.
A \textit{Lie color algebra} \(\mathfrak{g}\) is a \(\Gamma\)-graded \(\mathbb{R}\)-module \(\mathfrak{g} = \oplus_{\alpha \in \Gamma}\mathfrak{g}_\alpha\) with a bilinear map \([\, , ] : \mathfrak{g} \times \mathfrak{g} \longrightarrow \mathfrak{g}\) satisfying \([\mathfrak{g}_\alpha, \mathfrak{g}_\beta] \subseteq \mathfrak{g}_{\alpha+\beta}\) for every \(\alpha, \beta \in \Gamma\) and
\begin{equation}
[y, x] = -(\beta, \alpha)[x, y],
\end{equation}
\begin{equation}
(\gamma, \alpha)[x, [y, z]] + (\alpha, \beta)[y, [z, x]] + (\beta, \gamma)[z, [x, y]] = 0
\end{equation}
for homogeneous elements \(x \in \mathfrak{g}_\alpha\), \(y \in \mathfrak{g}_\beta\), \(z \in \mathfrak{g}_\gamma\) and \(\alpha, \beta, \gamma \in \Gamma\). (2.1) is called \textit{color skew-symmetry} and (2.2) is called \textit{the Jocobi color identity}. 


Let $A$ be a commutative-color algebra over $\mathbb{R}$ and $\delta\in\Gamma$. A $\mathbb{R}$-linear map $D:A\to A$ is called a \textit{homogeneous color derivation} of degree $\delta$ if it satisfies:
\[
D(a)\in A_{\delta+\alpha},
\]
\[
D(ab)=D(a)b+(\delta,\alpha)\,aD(b),
\]
for homogeneous elements $a\in A_{\alpha},\,b\in A_{\beta},\ \alpha,\beta\in\Gamma.$
Denote by $\mathrm{Der}^{\Gamma}_{\delta}(A)$ the $\mathbb{R}$-module spanned by all homogeneous color derivations of degree $\delta$ on $A$. \textit{The color derivation space} of $A$ is defined as
\[
\mathrm{Der}_{\mathbb{R}}^{\Gamma}(A):=\bigoplus_{\delta\in\Gamma}\mathrm{Der}^{\Gamma}_{\delta}(A).
\]
For any homogeneous color derivations $D\in\mathrm{Der}^{\Gamma}_{\delta}(A)$ and $D'\in\mathrm{Der}^{\Gamma}_{\delta'}(A)$, define \textit{the color bracket operation} on $\mathrm{Der}^{\Gamma}(A)$ by
\[
[D,D']:=DD'-(\delta,\delta')\,D'D.
\]
Then $\mathrm{Der}^{\Gamma}(A)$ equipped with the bracket operation forms a Lie color algebra \cite{ref30}.

\section{ Lie-Rinehart color algebra}
\begin{definition}
Let $A$ be a color-commutative algebra and $\mathfrak{g}=\bigoplus_{\alpha\in\Gamma}\mathfrak{g}_{\alpha}$ a Lie color algebra equipped with a (left) $A$-module structure. If the anchor map $\rho:\mathfrak{g}\to \mathrm{Der}^{\Gamma}_{\mathbb{R}}(A)$ is both an $A$-module homomorphism and a Lie color algebra homomorphism, satisfying
\[
[x,ay]=(\beta,\alpha)a[x,y]+\rho(x)(a)y,\quad \mbox{for all}\quad a\in A_{\alpha},\,x\in\mathfrak{g}_{\beta},\,y\in\mathfrak{g}_{\gamma},
\]
then the quintuple $(A,\mathfrak{g},\rho,\Gamma,(\cdot,\cdot))$ is called a Lie-Rinehart color algebra, which may be abbreviated as $(A,\mathfrak{g})$ when no confusion arises.
\end{definition}

If $\Gamma=\{0\}$ and \((0,0)=1\), then the Lie-Rinehart color algebra  is a Lie-Rinehart algebra. 

If $\Gamma=\mathbb{Z}_{2}$ and $(\beta,\gamma)=(-1)^{\beta\gamma}$, then the Lie-Rinehart color algebra is a Lie--Rinehart superalgebra. 

If $A=\mathbb{R}$, then the Lie-Rinehart color algebra is a Lie color algebra.

\noindent\textbf{Remark}: The definition differs in coefficients from the example given in \cite{ref5} when studying 3-Lie-Rinehart color algebras.


\begin{example}
Let $A$ be a color-commutative algebra and \(\mathfrak g\) be a \(\Gamma\)-graded abelian Lie color algebra equipped \ with a graded left \(A\)-module structure. Let $\rho:\mathfrak{g}\to \mathrm{Der}_{\mathbb{R}}^{\Gamma}(A)$ be an $A$-linear Lie color algebra homomorphism satisfying $\rho(x)=0$ for all $x\in\mathfrak{g}$. Then $(A,\mathfrak{g},\rho,\Gamma,(\cdot,\cdot))$ is a Lie-Rinehart color algebra.
\end{example}

\begin{example}
Let $A$ be a color-commutative algebra. Denote by $\operatorname{Der}_{\mathbb{R}}^{\Gamma}(A)$ the space of all $\Gamma$--graded color derivations of $A$, which is a Lie color algebra. Define a left $A$--module action by
\[
(a\cdot x)(b)=a\,x(b).
\]
Let \(\rho\) be the identity map on \(\operatorname{Der}_{\mathbb R}^{\Gamma}(A)\). Then the compatibility condition between the $A$-module structure on $\operatorname{Der}_{\mathbb{R}}^{\Gamma}(A)$ and the homomorphism of $A$-linear Lie color algebras holds. Hence, $(A,\operatorname{Der}_{\mathbb{R}}^{\Gamma}(A),\rho,\Gamma,(\cdot,\cdot))$ is a Lie-Rinehart color algebra.
\end{example}
\begin{definition}
Let $A$ be a color-commutative algebra. An $A$-module $M$ is a $\Gamma$-graded vector space $M=\bigoplus_{\gamma\in\Gamma}M_{\gamma}$ such that $A_{\alpha}\cdot M_{\beta}\subseteq M_{\alpha+\beta}$ for any $\alpha,\beta\in\Gamma$. A representation of a Lie-Rinehart color algebra $(A,\mathfrak{g},\rho,\Gamma,(\cdot,\cdot))$ consists of an $A$-module $M$ and an $A$-linear homomorphism $\phi:\mathfrak{g}\to\operatorname{End}^\Gamma_{\mathbb{R}}(M)$ of Lie color algebras, satisfying the identity:
\[
\phi(x)(am)=(\beta,\alpha)a\phi(x)(m)+\rho(x)(a)m
\]
for all $a\in A_\alpha$, $x\in\mathfrak{g}_\beta$, and $m\in M$.
The pair $(\phi,M)$ is called a Lie-Rinehart module. When no confusion arises, $M$ is called a Lie-Rinehart module.
\end{definition}
\begin{example}

 Let $\Gamma = \mathbb{Z}_2 \times \mathbb{Z}_2$. For $a = (a_1,a_2)$ and $b = (b_1,b_2)$ in $\mathbb{Z}_2 \times \mathbb{Z}_2$, set
$$ab := a_1b_2 - a_2b_1 \in \mathbb{Z}_2, (a,b): = (-1)^{ab}.$$
Let \(u\) and \(v\) be homogeneous generators of degrees
\((1,0)\) and \((0,1)\), respectively. We consider the algebra
$A:= \mathbb{R}\langle u,v \rangle / \langle uv + vu \rangle,$ where \(\mathbb{R}\langle u,v\rangle\)
is the free associative \(\mathbb{R}\)-algebra generated by \(u,v\),
and \(\langle u,v \rangle\) is the two-sided ideal generated by \(uv+vu\).
The algebra \(A\) has the $\mathbb{Z}_2 \times \mathbb{Z}_2$-grading
$A = \bigoplus\limits_{a\in\mathbb Z_2\times\mathbb Z_2} A^a,$
where
\[
A^a = \operatorname{span}_\mathbb{R}\{u^{2m+a_1}v^{2n+a_2} \mid m,n \ge 0\}.
\]
Moreover, $xy = (-1)^{ab}yx$ for $x \in A^a$ and $y \in A^b$. Then $A$ is a $\mathbb{Z}_2 \times \mathbb{Z}_2$-graded commutative algebra.
Set
\[\mathfrak{g} = A \otimes_\mathbb{R} \mathfrak{gl}(V), M = A \otimes_\mathbb{R} V,\]
with the total gradings
\[
\mathfrak{g}^a = \bigoplus_{b+c=a} A^b \otimes_\mathbb{R} \mathfrak{gl}(V)^c,\qquad
M^a = \bigoplus_{b+c=a} A^b \otimes_\mathbb{R} V^c.
\]
Define the left $A$-module structures by
\[
z(x \otimes f) = zx \otimes f,\qquad z(y \otimes w) = zy \otimes w
\]
for all $x,y,z \in A$, $f\in \mathfrak{gl}(V),w \in V$.
Define
\[
[x \otimes f, y \otimes g] = (-1)^{cb}xy \otimes [f,g],\qquad \rho = 0,
\]
for homogeneous elements
 $x \in A^a, y \in A^b,f \in \mathfrak{gl}(V)^c, g \in \mathfrak{gl}(V)^d.$
Then $(A,\mathfrak{g})$ is a Lie-Rinehart color algebra. Moreover, the map
\[
\phi \colon \mathfrak{g} \to \operatorname{End}_\mathbb{R}^\Gamma(M),\qquad
\phi(x \otimes f)(y \otimes w) = (-1)^{cb}xy \otimes f(w),
\]
 makes $(\phi,M)$ a Lie-Rinehart color module.
\end{example}

\begin{definition}
Let $(A,\mathfrak{g},\rho,\Gamma,(\cdot,\cdot))$ be a Lie-Rinehart color algebra. Suppose $\mathfrak{h}$ is a Lie color subalgebra of $\mathfrak{g}$ and also a graded $A$-submodule of $\mathfrak{g}$. Then $(A,\mathfrak{h},\rho,\Gamma,(\cdot,\cdot))$ is called a \textit{Lie-Rinehart color subalgebra} of $(A,\mathfrak{g},\rho,\Gamma,(\cdot,\cdot))$.
\end{definition}
\begin{example}
Let $A$ be a color-commutative algebra and \(\mathfrak g\) be a \(\Gamma\)-graded abelian Lie color algebra equipped \ with a graded left \(A\)-module structure. Let $\rho:\mathfrak{g}\to \mathrm{Der}_{\mathbb{R}}^{\Gamma}(A)$ be an $A$-linear Lie color algebra homomorphism satisfying $\rho(x)=0$ for all $x\in\mathfrak{g}$. Then $(A,\mathfrak{g},\rho,\Gamma,(\cdot,\cdot))$ is a Lie-Rinehart color algebra. Let $\mathfrak{h}$ be any graded $A$-submodule of $\mathfrak{g}$. Since $
[\mathfrak{h},\mathfrak{h}] = 0 \subseteq \mathfrak{h},$
$\mathfrak{h}$ is also a Lie color subalgebra of $\mathfrak{g}$. Hence  $(A,\mathfrak{h},\rho|_{\mathfrak{h}},\Gamma,(\cdot,\cdot))$
is a Lie-Rinehart color subalgebra of $(A,\mathfrak{g},\rho,\Gamma,(\cdot,\cdot))$.
\end{example}
\begin{definition}
 Let $(A,\mathfrak{g},\rho,\Gamma,(\cdot,\cdot))$ be a Lie-Rinehart color algebra. Suppose that $I$ is a Lie color ideal of $\mathfrak{g}$ and a graded $A$-submodule of $\mathfrak{g}$, satisfying $\rho(I)=0$. Then $I$ is called an ideal of $(A,\mathfrak{g},\rho,\Gamma,(\cdot,\cdot))$.
\end{definition}
\begin{example}
 Under the assumptions of Example 3.8, let $I$ be any $\Gamma$-graded $A$-submodule of $\mathfrak{g}$. Since
$[\mathfrak{g}, I] = 0 \subseteq I, \rho(I) = 0$, then
$I$ is an ideal of the Lie-Rinehart color algebra $(A, \mathfrak{g}, \rho, \Gamma, (\cdot,\cdot))$.
\end{example}
\begin{definition}
Let $(A,\mathfrak{g},\rho_1,\Gamma,(\cdot,\cdot))$ and $(A',\mathfrak{g'},\rho_2,\Gamma,(\cdot,\cdot))$ be Lie-Rinehart color algebras. A pair of morphisms $(\phi,f):\phi:A\to A'$ is a color algebra homomorphism, $f:\mathfrak{g}\to\mathfrak{g'}$ is a homomorphism of Lie color algebras and also a module homomorphism, satisfying
\[
\phi\circ\rho_1(x)(a)=\rho_2(f(x))(\phi(a)),
\]
for homogeneous elements a $\in A_\alpha, x\in\mathfrak{g}_\beta.$
Then the pair $(\phi,f)$ is called a homomorphism of Lie-Rinehart color algebras.
\end{definition}


 We now introduce  the category of Lie-Rinehart color algebras.

\begin{definition}
A category consists of objects and morphisms. The objects in the category of Lie-Rinehart color algebras are Lie-Rinehart color algebras, and the morphisms between objects are Lie-Rinehart color algebra homomorphisms $\varphi=(\phi,f)$.
\end{definition}

\begin{definition}
Let $(A,\mathfrak{g},\rho,\Gamma,(\cdot,\cdot))$ be a Lie-Rinehart color algebra and $\theta\in\Gamma$. A $\theta$-degree derivation of the Lie-Rinehart color algebra is a pair of homogeneous linear maps $(D,\delta)$ of degree $\theta$, where $D:\mathfrak{g}\to\mathfrak{g}$, $\delta:A\to A$, satisfying
\[
D(\mathfrak{g}_{\alpha})\subseteq \mathfrak{g}_{\alpha+\theta},\qquad \delta(A_{\beta})\subseteq A_{\beta+\theta},
\]
and the following identities hold:
\[
\delta(ab)=\delta(a)b+(\theta,\beta)a\delta(b);
\]
\[
D([x,y])=[D(x),y]+(\theta,\alpha)[x,D(y)];
\]
\[
D(ax)=\delta(a)x+(\theta,\beta)aD(x),
\]
for all homogeneous elements $a\in A_{\beta}$, $b\in A_{\gamma}$, $x\in\mathfrak{g}_{\alpha}$ and $y\in\mathfrak{g}_{\zeta}$.
The set of all $\theta$-degree Lie-Rinehart color derivations is denoted by $\operatorname{Der}_{\theta}^{\Gamma}(A,\mathfrak{g})$. Set
\[
\operatorname{Der}_{\mathbb{R}}^\Gamma(A,\mathfrak{g}):=\bigoplus_{\theta\in\Gamma}\operatorname{Der}^\Gamma_{\theta}(A,\mathfrak{g}).
\]
\end{definition}
For every homogeneous element $z \in \mathfrak{g}_\theta$, the pair $(\operatorname{ad}_z, \rho(z))$, where
\[
\operatorname{ad}_z(x) = [z, x], \quad \mbox{for all}\, x \in \mathfrak{g},
\]
is a derivation of degree $\theta$ of the Lie-Rinehart color algebra $(A, \mathfrak{g}, \rho, \Gamma, (\cdot,\cdot))$. Derivations of this form are called \textit{inner derivations}.
\section{The Universal Enveloping Algebras of Lie-Rinehart color algebras}
\subsection{Construction of the Universal Enveloping Algebra}
\begin{proposition}
Let $(A,\mathfrak{g},\rho,\Gamma,(\cdot,\cdot))$ be a Lie-Rinehart color algebra. Define the bracket operation on the graded $\mathbb{R}$-module $A\oplus\mathfrak{g}$ by
\[
[a+x,b+y]=[x,y]_{\mathfrak{g}}+\rho(x)(b)-(\alpha,\beta)\rho(y)(a),
\]
where $a+x\in(A\oplus\mathfrak{g})_{\alpha}, \ b+y\in(A\oplus\mathfrak{g})_{\beta}.$
Then \(A\oplus\mathfrak g\) is a Lie color algebra, denoted by $A\rtimes\mathfrak{g}$.
\end{proposition}
\begin{proof}
For $a+x\in(A\oplus\mathfrak g)_\alpha$ and
$b+y\in(A\oplus\mathfrak g)_\beta$, the defining formula gives
\[
[b+y,a+x]=-(\beta,\alpha)[a+x,b+y].
\]
The color Jacobi identity follows from the color Jacobi identity in
$\mathfrak g$ and the fact that $\rho$ is a homomorphism of Lie color
algebras.
\end{proof}

\noindent\textbf{Remark:} For \(a,a_1\in A\) and \(x\in\mathfrak g\), the left $A$-module action on $A\rtimes\mathfrak{g}$ is given by:
\[
a(a_1+x)=aa_1+ax.
\]

We next introduce the universal enveloping algebra of Lie-Rinehart color algebras. 
\begin{definition}
Let $(A,\mathfrak{g},\rho,\Gamma,(,))$ be a Lie-Rinehart color algebra. Its universal enveloping algebra is a triple $\big(U(A,\mathfrak{g}),\iota_A,\iota_{\mathfrak{g}}\big)$, where $U(A,\mathfrak{g})$ is a color associative algebra, $\iota_A\colon A\to U(A,\mathfrak{g})$ is a homomorphism of color algebras, and $\iota_{\mathfrak{g}}\colon \mathfrak{g}\to U(A,\mathfrak{g})^-$ is a homomorphism of Lie color algebras, satisfying the following compatibility conditions:
\[
\iota_A(a)\iota_{\mathfrak{g}}(x)=\iota_{\mathfrak{g}}(ax),\qquad
\iota_A(\rho(x)(a))=\iota_{\mathfrak{g}}(x)\iota_A(a)-(\beta,\alpha)\iota_A(a)\iota_{\mathfrak{g}}(x),
\]
for all homogeneous elements $a\in A_\alpha,\,x\in \mathfrak{g}_\beta$. Moreover, for any color algebra $B$, if there exists a homomorphism of color algebras $f\colon A\to B$ and a homomorphism of Lie color algebras $g\colon \mathfrak{g}\to B^-$,satisfing:
\[
g(ax)=f(a)g(x);
\]
\[
f(\rho(x)(a))=g(x)f(a)-(\beta,\alpha)f(a)g(x).
\]
for all $a\in A_{\alpha}$, $x\in\mathfrak{g}_{\beta}$,
then there exists a unique homomorphism of color algebras $\widetilde{F}:U(A,\mathfrak{g})\to B$ such that $
\widetilde{F}\circ\iota_{\mathfrak{g}}=g \quad \text{and} \quad \widetilde{F}\circ\iota_A=f.$
\end{definition}

We first prove the existence of this algebraic structure and verify its universal property.
\begin{theorem}
Every Lie-Rinehart color algebra admits a universal enveloping algebra.
\end{theorem}
\begin{proof}
Suppose $(A,\mathfrak{g},\rho,\Gamma,(\cdot,\cdot))$ is a Lie-Rinehart color algebra. Let \(A\rtimes\mathfrak g\) be the semidirect product Lie color algebra constructed above.
Denote by $U(A\rtimes\mathfrak{g})$ the universal enveloping algebra of $A\rtimes\mathfrak{g}$. Let $i:A\rtimes\mathfrak{g}\to U(A\rtimes\mathfrak{g})$ be the canonical map and $U^{+}(A\rtimes\mathfrak{g})$ be the associative algebra without unit generated by $i(A\rtimes\mathfrak{g})$. 
let \(I\) be the two-sided ideal of $U^{+}(A\rtimes\mathfrak{g})$  generated by
\[
i(a)i(a_1+x)-i(a(a_1+x)),
\]
where \(a,a_1\in A\) and \(x\in\mathfrak g\).

Set \(U(A,\mathfrak g)=U^{+}(A\rtimes\mathfrak g)/I\), and let
\(\pi:U^{+}(A\rtimes\mathfrak g)\to U(A,\mathfrak g)\)
be the quotient map. Define
\[
\iota_A:A\to U(A,\mathfrak{g}),\qquad
\iota_A(a)=\pi\bigl(i(a+0)\bigr),
\]
and
\[
\iota_{\mathfrak g}:\mathfrak g\to U(A,\mathfrak{g}),\qquad
\iota_{\mathfrak g}(x)=\pi\bigl(i(0+x)\bigr).
\]
Then we have
\[
\iota_A(a)\iota_{\mathfrak{g}}(x)=\iota_{\mathfrak{g}}(ax).
\]
\[
\iota_{\mathfrak{g}}(x)\iota_A(a)-(\beta,\alpha)\iota_A(a)\iota_{\mathfrak{g}}(x)=\iota_A(\rho(x)(a))
\]
Next we show that \((U(A,\mathfrak{g}),\iota_A,\iota_{\mathfrak g})\) is a universal enveloping algebra of \((A,\mathfrak g)\).
 Let $B$ be a color algebra. Suppose $f:A\to B$ is a homomorphism of color algebras, and $g:\mathfrak{g}\to B^-$ is a homomorphism of Lie color algebras, satisfing:
\[
g(ax)=f(a)g(x);
\]
\[
f(\rho(x)(a))=g(x)f(a)-(\beta,\alpha)f(a)g(x).
\]
for all $a\in A_{\alpha}$, $x\in\mathfrak{g}_{\beta}$.

 Define a linear map $
F:A\rtimes\mathfrak{g}\to B, F(a+z)=f(a)+g(z).$
It follows that $F$ is a homomorphism of
Lie color algebras.
By the universal property of the universal enveloping algebra $U(A\rtimes\mathfrak{g})$ of the Lie color algebra $A\rtimes\mathfrak{g}$, there exists a unique homomorphism of color algebras $\overline{F}:U(A\rtimes\mathfrak{g})\to B$ satisfying $\overline{F}\circ i=F$.
Since $I\subseteq\ker\overline{F}$, then $\overline{F}$ induces a homomorphism of color algebras $\widetilde{F}:U^{+}(A\rtimes\mathfrak{g})/I\to B$, such that $\widetilde{F}\circ\iota_{\mathfrak{g}}=g$ and $\widetilde{F}\circ\iota_A=f$. Since \(U(A,\mathfrak g)\) is generated by \(\iota_A(A)\) and \(\iota_{\mathfrak g}(\mathfrak g)\), then uniqueness of $\widetilde{F}$ follows.
\end{proof}
Following the theorem above, we obtain the following commutative diagram:
\[
\begin{tikzcd}
A \rtimes \mathfrak{g}
    \arrow[r, "i"] 
    \arrow[dr, "F"']
& U^+(A \rtimes \mathfrak{g}) 
    \arrow[r, "\pi"] 
    \arrow[d, "\bar{F}"']
& U(A,\mathfrak{g}) = U^+(A \rtimes \mathfrak{g})/I 
    \arrow[dl, "\widetilde{F}"] \\
& B
\end{tikzcd}
\]

\begin{corollary}
The map $\iota_A:A\to U(A,\mathfrak{g})$ is a monomorphism of color algebras.
\end{corollary}
\begin{proof}
 Let $B=\operatorname{End}_{\mathbb{R}}^{\Gamma}(A)$. Define the map $f:A\to\operatorname{End}_{\mathbb{R}}^{\Gamma}(A)$ by $f(a)(b)=a\cdot b$ for all $a,b\in A$. Then $f$ is a homomorphism of color algebras. Define the map $g:\mathfrak{g}\to\mathfrak{gl}_{\mathbb{R}}^{\Gamma}(A)$ by $g(x)(a):=\rho(x)(a)$, for all $x\in\mathfrak{g},a\in A$. Then $g$ is a homomorphism of Lie color algebras. The following relations hold:
\[
g(ax)=f(a)g(x);
\]
\[
f(\rho(x)(a))=g(x)f(a)-(\beta,\alpha)f(a)g(x)
\]
for all $a\in A_{\alpha},x\in\mathfrak{g}_{\beta}$.
Then there exists a unique homomorphism of color algebras $\widetilde{F}:U(A,\mathfrak{g})\to\operatorname{End}_{\mathbb{R}}^{\Gamma}(A)$ such that $\widetilde{F}\circ\iota_A=f$. Since $f$ is a monomorphism, then we conclude that $\iota_A$ is a monomorphism.
\end{proof}
\subsection{Properties of the Universal Enveloping Algebra}
\begin{theorem}
Let $(U(A,\mathfrak{g}),\iota_A,\iota_{\mathfrak{g}})$ and $(V(A,\mathfrak{g}),j_A,j_{\mathfrak{g}})$ be two universal enveloping algebras of $(A,\mathfrak{g})$. Then there exists a unique isomorphism of color algebras $h:U(A,\mathfrak{g})\to V(A,\mathfrak{g})$ satisfying $h\circ\iota_A=j_A$ and $h\circ\iota_{\mathfrak{g}}=j_{\mathfrak{g}}$.
\end{theorem}
\begin{proof}
 Since $(U(A,\mathfrak{g}),\iota_A,\iota_{\mathfrak{g}})$ is the universal enveloping algebra of $(A,\mathfrak{g})$, and the color algebra $(V(A,\mathfrak{g}),j_A,j_{\mathfrak{g}})$ satisfies,
\[
j_A(a)j_{\mathfrak{g}}(x)=j_{\mathfrak{g}}(ax);
\]
\[
j_A(\rho(x)(a))=[j_{\mathfrak{g}}(x),j_A(a)],
\]
for all homogeneous elements $a\in A$, $x\in\mathfrak{g}$, 
then there exists a unique homomorphism of color algebras $\widetilde{F}:U(A,\mathfrak{g})\to V(A,\mathfrak{g})$ such that $\widetilde{F}\circ\iota_{\mathfrak{g}}=j_{\mathfrak{g}}$ and $\widetilde{F}\circ\iota_A=j_A$.

Similarly, since $(V(A,\mathfrak{g}),j_A,j_{\mathfrak{g}})$ is the universal enveloping algebra of $(A,\mathfrak{g})$, then there exists a unique homomorphism of color algebras $\widetilde{G}:V(A,\mathfrak{g})\to U(A,\mathfrak{g})$ such that $\widetilde{G}\circ j_A=\iota_A$ and $\widetilde{G}\circ j_{\mathfrak{g}}=\iota_{\mathfrak{g}}$.
Then we have $\widetilde{F}\circ\iota_{\mathfrak{g}}=j_{\mathfrak{g}}$ and $\widetilde{G}\circ j_{\mathfrak{g}}=\iota_{\mathfrak{g}}$. 
Furthermore,
\[
\begin{aligned}
(\widetilde F\circ\widetilde G)\circ j_A &= j_A,
&\qquad
(\widetilde F\circ\widetilde G)\circ j_{\mathfrak g} &= j_{\mathfrak g},\\
(\widetilde G\circ\widetilde F)\circ \iota_A &= \iota_A,
&\qquad
(\widetilde G\circ\widetilde F)\circ \iota_{\mathfrak g} &= \iota_{\mathfrak g}.
\end{aligned}
\]
 By the uniqueness property of universal enveloping algebras, we get $\widetilde{F}\circ\widetilde{G}=\operatorname{id}_{V(A,\mathfrak{g})}$. Similarly, $\widetilde{G}\circ\widetilde{F}=\operatorname{id}_{U(A,\mathfrak{g})}$. Hence $\widetilde{F}$ and $\widetilde{G}$ are mutual inverses, so $\widetilde{G}=\widetilde{F}^{-1}$.
  Let $h=\widetilde{F}$. Then the conclusion holds.
\end{proof}
\begin{theorem}
Let $(U(A,\mathfrak{g}),\iota_A,\iota_{\mathfrak{g}})$ be the universal enveloping algebra of $(A,\mathfrak{g})$. Then $U(A,\mathfrak{g})=\operatorname{alg}_A\big(\iota_{\mathfrak{g}}(\mathfrak{g})\big)$.
\end{theorem}
\begin{proof}
Let $C=\operatorname{alg}_R\big(\iota_A(A)\cup\iota_{\mathfrak{g}}(\mathfrak{g})\big)$ be the unital $\mathbb{R}$-color subalgebra of $U(A,\mathfrak{g})$ generated by $\iota_A(A)$ and $\iota_{\mathfrak{g}}(\mathfrak{g})$. 
Let $j:C\to U(A,\mathfrak{g})$ denote the inclusion map.

Define maps $f_A:A\to C$ by $f_A(a)=\iota_A(a)$ and $f_{\mathfrak{g}}:\mathfrak{g}\to C$ by $f_{\mathfrak{g}}(x)=\iota_{\mathfrak{g}}(x)$. Then $f_A$ is a homomorphism of color algebras and $f_{\mathfrak{g}}$ is a homomorphism of Lie color algebras. Then we have
\[
f_{\mathfrak{g}}(ax)=\iota_{\mathfrak{g}}(ax)=\iota_A(a)\iota_{\mathfrak{g}}(x)=f_A(a)f_{\mathfrak{g}}(x).
\]
\[
\big[f_{\mathfrak{g}}(x),f_A(a)\big]=\big[\iota_{\mathfrak{g}}(x),\iota_A(a)\big]=\iota_A(\rho(x)(a))=f_A(\rho(x)(a)).
\]
 for all homogeneous elements $a\in A_{\alpha}$ and $x\in\mathfrak{g}_{\beta}$.
By the universal property of $U(A,\mathfrak{g})$, there exists a unique homomorphism of color algebras $F:U(A,\mathfrak{g})\to C$ such that $F\circ\iota_A=f_A$ and $F\circ\iota_{\mathfrak{g}}=f_{\mathfrak{g}}$. Now consider the composite map $j\circ F:U(A,\mathfrak{g})\to U(A,\mathfrak{g})$. 
\[
((j\circ F)\circ\iota_A)(a)=j\big(f_A(a)\big)=\iota_A(a).
\]
\[
((j\circ F)\circ\iota_{\mathfrak{g}})(x)=j\big(f_{\mathfrak{g}}(x)\big)=\iota_{\mathfrak{g}}(x),
\]
for all $a\in A$, $x\in\mathfrak{g}$. The identity map satisfies $\operatorname{id}_{U(A,\mathfrak{g})}\circ\iota_A=\iota_A$ and $\operatorname{id}_{U(A,\mathfrak{g})}\circ\iota_{\mathfrak{g}}=\iota_{\mathfrak{g}}$. By universality, $j\circ F=\operatorname{id}_{U(A,\mathfrak{g})}$. It follows that
\[
U(A,\mathfrak{g})=j\big(F(U(A,\mathfrak{g}))\big)=j(C)=C.
\]
Thus
\[
U(A,\mathfrak{g})=\operatorname{alg}_R\big(\iota_A(A)\cup\iota_{\mathfrak{g}}(\mathfrak{g})\big).
\]
Since $\iota_A$ is injective, then we may identify $A$ with $\iota_A(A)$, which further yields $U(A,\mathfrak{g})=\operatorname{alg}_A\big(\iota_{\mathfrak{g}}(\mathfrak{g})\big)$.
\end{proof}
\begin{theorem}
Let $(A,\mathfrak{g},\rho_{\mathfrak{g}},\Gamma,(\cdot,\cdot))$ and $(B,\mathfrak{h},\rho_{\mathfrak{h}},\Gamma,(\cdot,\cdot))$ be two Lie-Rinehart color algebras. Let the pair $(\phi,f)\colon (A,\mathfrak{g},\rho_{\mathfrak{g}},\Gamma,(\cdot,\cdot))\to (B,\mathfrak{h},\rho_{\mathfrak{h}},\Gamma,(\cdot,\cdot))$ be a homomorphism of Lie-Rinehart color algebras. Let $(U(A,\mathfrak{g}),\iota_A,\iota_{\mathfrak{g}})$ and $(U(B,\mathfrak{h}),\iota_B,\iota_{\mathfrak{h}})$ be the universal enveloping algebras of $(A,\mathfrak{g})$ and $(B,\mathfrak{h})$, respectively. Then there exists a unique homomorphism of color algebras $U(\phi,f)\colon U(A,\mathfrak{g})\to U(B,\mathfrak{h})$ making the following diagrams commute.
\[
\newcommand{\arr}[1]{\xrightarrow{\makebox[2.5em]{$#1$}}}
\begin{array}{c@{\qquad}c}
\begin{array}{ccc}
A & \arr{\phi} & B \\
\downarrow^{\iota_A} & & \downarrow^{\iota_B}\\
U(A,\mathfrak{g}) & \arr{U(\phi,f)} & U(B,\mathfrak{h})
\end{array}
&
\begin{array}{ccc}
\mathfrak{g} & \arr{f} & \mathfrak{h} \\
\downarrow^{\iota_{\mathfrak{g}}} & & \downarrow^{\iota_{\mathfrak{h}}}\\
U(A,\mathfrak{g}) & \arr{U(\phi,f)} & U(B,\mathfrak{h})
\end{array}
\end{array}
\]
That is,
\[
U(\phi,f)\circ\iota_A=\iota_B\circ\phi,\qquad U(\phi,f)\circ\iota_{\mathfrak{g}}=\iota_{\mathfrak{h}}\circ f.
\]
\end{theorem}
\begin{proof}
The map $\iota_B\circ\phi\colon A\to U(B,\mathfrak{h})$ is a homomorphism of color algebras, and $\iota_{\mathfrak{h}}\circ f\colon \mathfrak{g}\to U(B,\mathfrak{h})^-$ is a homomorphism of Lie color algebras.
Set $F_A:=\iota_B\circ\phi$ and $F_{\mathfrak{g}}:=\iota_{\mathfrak{h}}\circ f$. Then we have
\[
F_{\mathfrak{g}}(ax)=\iota_{\mathfrak{h}}(f(ax))=\iota_{\mathfrak{h}}(\phi(a)f(x))=\iota_B(\phi(a))\iota_{\mathfrak{h}}(f(x))=F_A(a)F_{\mathfrak{g}}(x),
\]
\[
\big[F_{\mathfrak{g}}(x),F_A(a)\big]=\big[\iota_{\mathfrak{h}}(f(x)),\iota_B(\phi(a))\big]
=\iota_B\big(\rho_{\mathfrak{h}}(f(x))(\phi(a))\big)
=F_A\big(\rho_{\mathfrak{g}}(x)(a)\big),
\]
 for any homogeneous elements $a\in A_{\alpha},x\in\mathfrak{g}_{\beta}$.
Therefore, by the universal property of $U(A,\mathfrak{g})$, there exists a unique homomorphism of color algebras $$U(\phi,f)\colon U(A,\mathfrak{g})\to U(B,\mathfrak{h})$$ such that $U(\phi,f)\circ\iota_A=\iota_B\circ\phi$ and $U(\phi,f)\circ\iota_{\mathfrak{g}}=\iota_{\mathfrak{h}}\circ f$.
  
\end{proof}
\noindent\textbf{Remark:} $U(\operatorname{id}_A,\operatorname{id}_{\mathfrak{g}})=\operatorname{id}_{U(A,\mathfrak{g})}$. 
\begin{corollary}
Let $(B,\mathfrak{h},\rho_{\mathfrak{h}},\Gamma,(\cdot,\cdot))$ and
$(C,\mathfrak{l},\rho_{\mathfrak{l}},\Gamma,(\cdot,\cdot))$ be two Lie-Rinehart color algebras. 
Let $(U(C,\mathfrak{l}),\iota_C,\iota_{\mathfrak{l}})$ be the universal enveloping algebra of $(C,\mathfrak{l})$. If $(\phi,f)\colon(A,\mathfrak{g})\to(B,\mathfrak{h})$ and $(\psi,k)\colon(B,\mathfrak{h})\to(C,\mathfrak{l})$ are homomorphisms of Lie-Rinehart color algebras, then
\[
U(\psi,k)\circ U(\phi,f)=U(\psi\circ\phi,k\circ f).
\]
\end{corollary}
\begin{proof}
 Consider the composite morphism $$U(\psi,k)\circ U(\phi,f)\colon U(A,\mathfrak{g})\to U(C,\mathfrak{l}).$$ Then we obtain the following commutative diagrams.
 \[
\begin{array}{ccccc}
U(A,\mathfrak g)
& \xrightarrow{\makebox[1.5cm]{$\scriptstyle U(\phi,f)$}}&
U(B,\mathfrak h)
& \xrightarrow{\makebox[1.5cm]{$\scriptstyle U(\psi,k)$}} &
U(C,\mathfrak l)
\\[6pt]
{\scriptstyle\iota_\mathfrak g}\uparrow
&&{\scriptstyle\iota_\mathfrak h}\uparrow
&&{\scriptstyle\iota_\mathfrak l}\uparrow
\\[6pt]
\mathfrak g
& \xrightarrow{\makebox[1.5cm]{$\scriptstyle f$}} &
\mathfrak h
& \xrightarrow{\makebox[1.5cm]{$\scriptstyle k$}} &
\mathfrak l
\end{array}
\]
and
\[
\begin{array}{ccccc}
U(A,\mathfrak g)
& \xrightarrow{\makebox[1.5cm]{$\scriptstyle U(\phi,f)$}} &
U(B,\mathfrak h)
& \xrightarrow{\makebox[1.5cm]{$\scriptstyle U(\psi,k)$}} &
U(C,\mathfrak l)
\\[6pt]
{\scriptstyle\iota_A}\uparrow
&&{\scriptstyle\iota_B}\uparrow
&&{\scriptstyle\iota_C}\uparrow
\\[6pt]
A
& \xrightarrow{\makebox[1.5cm]{$\scriptstyle\phi$}} &
B
& \xrightarrow{\makebox[1.5cm]{$\scriptstyle\psi$}} &
C
\end{array}
\]
Furthermore, we have
\[
\big(U(\psi,k)\circ U(\phi,f)\big)\circ\iota_A
=U(\psi,k)\circ\iota_B\circ\phi
=\iota_C\circ\psi\circ\phi.
\]
\[
\big(U(\psi,k)\circ U(\phi,f)\big)\circ\iota_{\mathfrak{g}}
=U(\psi,k)\circ\iota_{\mathfrak{h}}\circ f
=\iota_{\mathfrak{l}}\circ k\circ f.
\]
Hence, $U(\psi,k)\circ U(\phi,f)=U(\psi\circ\phi,k\circ f)$.
\end{proof}
\begin{theorem}
Let $(A,\mathfrak{g},\rho,\Gamma,(\cdot,\cdot))$ be a Lie-Rinehart color algebra and $J$ be a Lie-Rinehart color ideal of $(A,\mathfrak{g})$. Let $\iota_A:A\to U(A,\mathfrak{g})$ and $\iota_{\mathfrak{g}}:\mathfrak{g}\to U(A,\mathfrak{g})^-$ be the canonical maps. Let $I$ be the two-sided ideal of $U(A,\mathfrak{g})$ generated by $\iota_{\mathfrak{g}}(J)$. Then $(U(A,\mathfrak{g})/I,\bar{\iota}_A,j)$ is the universal enveloping algebra of the Lie-Rinehart color algebra $(A,\mathfrak{g}/J,\bar{\rho},\Gamma,(\cdot,\cdot))$, where $\bar{\iota}_A$ and $j$ are the corresponding induced maps.
\end{theorem}
\begin{proof}
 Denote by $q:\mathfrak{g}\to\mathfrak{g}/J$ and $\pi:U(A,\mathfrak{g})\to U(A,\mathfrak{g})/I$ the canonical projections. 
 We define $a(x+J)=ax+J$ for all $a\in A_{\alpha},x\in\mathfrak{g}_{\beta}$. Then $\mathfrak{g}/J$ is an $A$-module. Define the map $\bar{\rho}:\mathfrak{g}/J\to\operatorname{Der}(A)$ by $\bar{\rho}(x+J)=\rho(x)$. Then we have
\[
[x+J,a(y+J)]=(\beta,\alpha)a\big([x+J,y+J]\big)+\bar{\rho}(x+J)(a)(y+J)
\]
for all $x+J\in(\mathfrak{g}/J)_{\beta}$, $y+J\in(\mathfrak{g}/J)_{\gamma}$, $a\in A_{\alpha}$. Hence $(A,\mathfrak{g}/J,\bar{\rho},\Gamma,(\cdot,\cdot))$ is a Lie--Rinehart color algebra.

Let $\pi\circ\iota_{\mathfrak{g}}:\mathfrak{g}\to \big(U(A,\mathfrak{g})/I\big)^-$. Since 
$J\subseteq\ker(\pi\circ\iota_{\mathfrak{g}})$, then there exists a unique homomorphism of Lie color algebras $j:\mathfrak{g}/J\to \big(U(A,\mathfrak{g})/I\big)^-$ such that $j\circ q=\pi\circ\iota_{\mathfrak{g}}$. Therefore, 
\begin{align*}
j(x+J)&=j(q(x))=(j\circ q)(x)=(\pi\circ\iota_{\mathfrak{g}})(x)\\
&=\pi(\iota_{\mathfrak{g}}(x))=\iota_{\mathfrak{g}}(x)+I.
\end{align*}
for all elements $x\in\mathfrak{g}_{\beta}$.
Set $$\bar{\iota}_A=\pi\circ\iota_A:A\to U(A,\mathfrak{g})\to U(A,\mathfrak{g})/I$$ with $\bar{\iota}_A(a)=\iota_A(a)+I$ for all $a\in A_{\alpha}$. Then we have
\begin{align*}
j(a(x+J))&=j(ax+J)=j(q(ax))=\pi(\iota_{\mathfrak{g}}(ax))=\pi(\iota_A(a)\iota_{\mathfrak{g}}(x))\\
&=\pi(\iota_A(a))\pi(\iota_{\mathfrak{g}}(x))=\bar{\iota}_A(a)j(x+J).
\end{align*}
\begin{align*}
\bar{\iota}_A\big(\bar{\rho}(x+J)(a)\big)&=\pi\big(\iota_A(\rho(x)(a))\big)=\pi(\iota_{\mathfrak{g}}(x))\pi(\iota_A(a))-(\beta,\alpha)\pi(\iota_A(a))\pi(\iota_{\mathfrak{g}}(x))\\
&=j(x+J)\bar{\iota}_A(a)-(\beta,\alpha)\bar{\iota}_A(a)j(x+J).
\end{align*}

We now verify the universal property. Let $B$ be a color algebra. Let $f:A\to B$ be a homomorphism of color algebras and $g:\mathfrak{g}/J\to B^-$ a homomorphism of Lie color algebras. Assume that 
\[
g\big(a(x+J)\big)=f(a)g(x+J),
\]
\[
f\big(\bar{\rho}(x+J)(a)\big)=g(x+J)f(a)-(\beta,\alpha)f(a)g(x+J).
\]
for any homogeneous elements $a\in A_{\alpha},x\in\mathfrak{g}_{\beta}$.
The composite map $g\circ q:\mathfrak{g}\to B^-$ is a homomorphism of Lie color algebras. 
\[
(g\circ q)(ax)=g(q(ax))=g(ax+J)=g\big(a(x+J)\big)=f(a)g(x+J)=f(a)(g\circ q)(x).
\]
\[
f(\rho(x)(a))=f\big(\bar{\rho}(q(x))(a)\big)=(g\circ q)(x)f(a)-(\beta,\alpha)f(a)(g\circ q)(x).
\]
By the universal property of $U(A,\mathfrak{g})$, there exists a unique homomorphism of color algebras $F:U(A,\mathfrak{g})\to B$ such that
\[
F\circ\iota_A=f,\qquad F\circ\iota_{\mathfrak{g}}=g\circ q.
\]
Hence $\iota_{\mathfrak{g}}(J)\subseteq\ker F$. Since $I$ is the two-sided ideal generated by $\iota_{\mathfrak{g}}(J)$, then we get $I\subseteq\ker F$. Therefore there exists a unique homomorphism of color algebras $\widetilde{F}:U(A,\mathfrak{g})/I\to B$ satisfying $\widetilde{F}\circ\pi=F$. Then
\[
\widetilde{F}\circ\bar{\iota}_A=\widetilde{F}\circ\pi\circ\iota_A=F\circ\iota_A=f.
\]
\[
\widetilde{F}\circ j\circ q=\widetilde{F}\circ\pi\circ\iota_{\mathfrak{g}}=F\circ\iota_{\mathfrak{g}}=g\circ q.
\]
Since $q$ is surjective, then $\widetilde{F}\circ j=g$. The uniqueness of $\widetilde{F}$ is obvious.
There exists a unique homomorphism of color algebras $\widetilde{F}:U(A,\mathfrak{g})/I\to B$ such that
\[
\widetilde{F}\circ\bar{\iota}_A=f \quad\text{and}\quad \widetilde{F}\circ j=g.
\]
Thus $(U(A,\mathfrak{g})/I,\bar{\iota}_A,j)$ is the universal enveloping algebra of the Lie--Rinehart color algebra $(A,\mathfrak{g}/J)$, and there is color algebra isomorphism
\[
U(A,\mathfrak{g}/J)\cong U(A,\mathfrak{g})/I.
\]
\end{proof}
\begin{theorem}
Let $(A,\mathfrak{g},\Gamma,(\cdot,\cdot),\rho)$ be a Lie-Rinehart color algebra, and suppose that $\mathfrak{g}$ is a faithful $A$-module. Let $(D,\delta)\in\operatorname{Der}_{\theta}^{\Gamma}(A,\mathfrak{g})$ be a $\theta$-degree Lie-Rinehart color derivation, where $\theta\in\Gamma$. Then there exists a unique $\theta$-degree derivation $\overline{D}:U(A,\mathfrak{g})\to U(A,\mathfrak{g})$ such that $\overline{D}\circ\iota_A=\iota_A\circ\delta$ and $\overline{D}\circ\iota_{\mathfrak{g}}=\iota_{\mathfrak{g}}\circ D$. Namely, the following commutative diagram commutes.
\[
\begin{array}{ccc}
A & \xrightarrow{\ \delta\ } & A\\[6pt]
{\scriptstyle\iota_A}\downarrow
&&\downarrow{\scriptstyle\iota_A}\\[6pt]
U(A,\mathfrak g)
&\xrightarrow{\ \overline D\ }&
U(A,\mathfrak g)
\end{array}
\qquad
\begin{array}{ccc}
\mathfrak g & \xrightarrow{\ D\ } & \mathfrak g\\[6pt]
{\scriptstyle\iota_{\mathfrak g}}\downarrow
&&\downarrow{\scriptstyle\iota_{\mathfrak g}}\\[6pt]
U(A,\mathfrak g)
&\xrightarrow{\ \overline D\ }&
U(A,\mathfrak g)
\end{array}
\]
\end{theorem}
\begin{proof}
We expand $D([x,ay])$ in different ways and compare the resulting
expressions to obtain
\[
\bigl(\delta(\rho(x)(a))-\rho(D(x))(a)
      -(\theta,\alpha)\rho(x)(\delta(a))\bigr)y=0.
\]
for all homogeneous elements $x\in\mathfrak{g}_{\alpha}$, $a\in A_{\beta}$, $y\in\mathfrak{g}_{\gamma},$ where $\alpha,\beta,\gamma\in\Gamma$ are arbitrary. 
Since $\mathfrak{g}$ is a faithful $A$-module, we have
\[
\delta(\rho(x)(a))=\rho(D(x))(a)+(\theta,\alpha)\rho(x)(\delta(a)).
\]

Define a linear map $\sigma_{\theta}:U(A,\mathfrak{g})\to U(A,\mathfrak{g})$ by $\sigma_{\theta}(u)=(\theta,\gamma)u$, for $u\in U(A,\mathfrak{g})_{\gamma}$. Then $\sigma_{\theta}$ is a color algebra automorphism of $U(A,\mathfrak{g})$.

Define the algebra $$M=\left\{\begin{pmatrix}\sigma_{\theta}(u)&v\\0&u\end{pmatrix}\,\bigg|\,u,v\in U\right\}.$$ Then $M=\bigoplus_{\gamma\in\Gamma}M_{\gamma}$ is a color algebra, where $$M_{\gamma}=\left\{\begin{pmatrix}\sigma_{\theta}(u)&v\\0&u\end{pmatrix}\,\bigg|\,u\in U_{\gamma},\,v\in U_{\gamma+\theta}\right\}.$$
Define maps $F_A:A\to M$ and $F_{\mathfrak{g}}:\mathfrak{g}\to M$ as follows:
\[
F_A(a)=\begin{pmatrix}\sigma_{\theta}(\iota_A(a))&\iota_A(\delta(a))\\0&\iota_A(a)\end{pmatrix},
\qquad
F_{\mathfrak{g}}(x)=\begin{pmatrix}\sigma_{\theta}(\iota_{\mathfrak{g}}(x))&\iota_{\mathfrak{g}}(D(x))\\0&\iota_{\mathfrak{g}}(x)\end{pmatrix}.
\]
Then we obtain that
$F_A$ is a homomorphism of color algebras and $F_{\mathfrak{g}}$ is a homomorphism of Lie color algebras.


Let $a\in A_{\beta}$ and $x\in\mathfrak{g}_{\alpha}$ be homogeneous elements. By direct computations, we have
\[
F_{\mathfrak{g}}(ax)=F_A(a)F_{\mathfrak{g}}(x).
\]
\[
F_{\mathfrak{g}}(x)F_A(a)-(\alpha,\beta)F_A(a)F_{\mathfrak{g}}(x)=F_A(\rho(x)(a)).
\]
By the universal property of $U(A,\mathfrak{g})$, there exists a unique homomorphism of color algebras $\widetilde{F}:U(A,\mathfrak{g})\to M$ such that $\widetilde{F}\circ\iota_A=F_A$ and $\widetilde{F}\circ\iota_{\mathfrak{g}}=F_{\mathfrak{g}}$.

Define the map $$p:M\to U(A,\mathfrak{g}), p\left(\begin{pmatrix}\sigma_{\theta}(u)&v\\0&u\end{pmatrix}\right)=u.$$ Clearly $p$ is a homomorphism of graded associative algebras.
\[
p\circ\widetilde{F}\circ\iota_A=p\circ F_A=\iota_A,\qquad p\circ\widetilde{F}\circ\iota_{\mathfrak{g}}=p\circ F_{\mathfrak{g}}=\iota_{\mathfrak{g}}.
\]
The identity map satisfies $\operatorname{id}_U\circ\iota_A=\iota_A$ and $\operatorname{id}_U\circ\iota_{\mathfrak{g}}=\iota_{\mathfrak{g}}$. Therefore $p\circ\widetilde{F}=\operatorname{id}_U$.

Define the map $$q:M\to U(A,\mathfrak{g}), q\left(\begin{pmatrix}\sigma_{\theta}(u)&v\\0&u\end{pmatrix}\right)=v.$$ Set $\overline{D}=q\circ\widetilde{F}: U(A,\mathfrak{g})\to  U(A,\mathfrak{g})$.
For any $u\in U(A,\mathfrak{g})$, since $\widetilde{F}(u)\in M$, there exist unique $u'\in U(A,\mathfrak{g})$ and $v\in U$ such that
\[
\widetilde{F}(u)=\begin{pmatrix}\sigma_{\theta}(u')&v\\0&u'\end{pmatrix}.
\]
where $u'=p(\widetilde{F}(u))=u$, $v=q(\widetilde{F}(u))=\overline{D}(u)$. Hence
\[
\widetilde{F}(u)=\begin{pmatrix}\sigma_{\theta}(u)&\overline{D}(u)\\0&u\end{pmatrix}.
\]
Thus
\[
\overline{D}\circ\iota_A=q\circ\widetilde{F}\circ\iota_A=q\circ F_A=\iota_A\circ\delta.
\]
\[
\overline{D}\circ\iota_{\mathfrak{g}}=q\circ\widetilde{F}\circ\iota_{\mathfrak{g}}=q\circ F_{\mathfrak{g}}=\iota_{\mathfrak{g}}\circ D.
\]
Take homogeneous elements $u\in U(A,\mathfrak{g})_{\alpha}$, $v\in U(A,\mathfrak{g})_{\beta}$. Since $\widetilde{F}$ is a homomorphism of associative algebras, we have $$\widetilde{F}(uv)=\widetilde{F}(u)\widetilde{F}(v).$$ Therefore
\[
\begin{pmatrix}\sigma_{\theta}(uv)&\overline{D}(uv)\\0&uv\end{pmatrix}
=\begin{pmatrix}\sigma_{\theta}(u)&\overline{D}(u)\\0&u\end{pmatrix}
\begin{pmatrix}\sigma_{\theta}(v)&\overline{D}(v)\\0&v\end{pmatrix}
=\begin{pmatrix}\sigma_{\theta}(uv)&\sigma_{\theta}(u)\overline{D}(v)+\overline{D}(u)v\\0&uv\end{pmatrix}.
\]
That is,
\[
\overline{D}(uv)=\sigma_{\theta}(u)\overline{D}(v)+\overline{D}(u)v=\overline{D}(u)v+(\theta,\alpha)u\overline{D}(v).
\]
Hence $\overline{D}$ is a $\theta$-degree derivation of $U(A,\mathfrak{g})$. Finally, since $U(A,\mathfrak{g})$ is generated by $\iota_A(A)\cup\iota_{\mathfrak{g}}(\mathfrak{g})$, and the values of $\overline{D}$ on these generators are uniquely determined, $\overline{D}$ is unique.

There exists a unique $\theta$-degree associative algebra derivation $$\overline{D}:U(A,\mathfrak{g})\to U(A,\mathfrak{g})$$ satisfying $\overline{D}\circ\iota_A=\iota_A\circ\delta$ and $\overline{D}\circ\iota_{\mathfrak{g}}=\iota_{\mathfrak{g}}\circ D$.
\end{proof}
We next turn to modules over Lie-Rinehart color algebras and show that their actions extend uniquely to actions of the corresponding universal enveloping algebras.
\begin{theorem}
Let $(U(A,\mathfrak{g}),\iota_A,\iota_{\mathfrak{g}})$ be the universal enveloping algebra of $(A,\mathfrak{g})$. Let $(\phi,M)$ be a Lie-Rinehart color module. Then there exists a unique graded $U(A,\mathfrak{g})$-module structure on $M$ such that
\[
\iota_A(a)\cdot m = am,\qquad \iota_{\mathfrak{g}}(x)\cdot m=\phi(x)(m)
\]
for all $a\in A$, $x\in\mathfrak{g}$, $m\in M$. Furthermore, suppose $N\subseteq M$ is a graded $A$-submodule. Then $N$ is a Lie-Rinehart color submodule if and only if $N$ is a $U(A,\mathfrak{g})$-submodule of $M$.
\end{theorem}
\begin{proof}
 Define $\lambda:A\to\operatorname{End}^\Gamma_\mathbb{R}(M)$ by $\lambda(a)(m)=am$ for all $a\in A$ and $m\in M$. Then $\lambda$ is a graded associative algebra homomorphism, and $\phi:\mathfrak{g}\to\operatorname{gl}^\Gamma_\mathbb{R}(M)$ is a homomorphism of Lie color algebras.

For all homogeneous elements $a\in A_{\alpha}$, $x\in\mathfrak{g}_{\beta}$ and $m\in M$, we have
\[
\phi(ax)(m)=a\phi(x)(m)=\lambda(a)\phi(x)(m).
\]
\[
\big(\phi(x)\lambda(a)-(\beta,\alpha)\lambda(a)\phi(x)\big)(m)=\phi(x)(am)-(\beta,\alpha)a\phi(x)(m)=\lambda(\rho(x)(a))(m).
\]
Therefore, the pair $(\lambda,\phi)$ satisfies the universal property of the universal enveloping algebra $U(A,\mathfrak{g})$. Hence there exists a unique graded associative algebra homomorphism $\overline{\phi}:U(A,\mathfrak{g})\to\operatorname{End}^\Gamma_\mathbb{R}(M)$  such that $\overline{\phi}\circ\iota_A=\lambda$ and $\overline{\phi}\circ\iota_{\mathfrak{g}}=\phi$.

Define $u\cdot m=\overline{\phi}(u)(m)$ for all $u\in U(A,\mathfrak{g})$, $m\in M$. For all $a\in A_{\alpha}$, $x\in\mathfrak{g}_{\beta}$ and $m\in M$,
\[
\iota_A(a)\cdot m=\overline{\phi}(\iota_A(a))(m)=\lambda(a)(m)=am,
\]
\[
\iota_{\mathfrak{g}}(x)\cdot m=\overline{\phi}(\iota_{\mathfrak{g}}(x))(m)=\phi(x)(m).
\]
By the uniqueness of $\overline{\phi}$ in the universal property, the above $U(A,\mathfrak{g})$-module structure is also uniquely determined.

Let $N\subseteq M$ be an $A$-submodule. If $N$ is a $U(A,\mathfrak{g})$-submodule, then we have $\phi(x)(n)=\iota_{\mathfrak{g}}(x)\cdot n\in N$, for all $x\in\mathfrak{g}$ and $n\in N$. Thus $\phi(\mathfrak{g})(N)\subseteq N$. Since $N$ is a graded $A$-submodule, $N$ is a Lie--Rinehart color submodule.

Conversely, assume $N$ is a Lie--Rinehart color submodule. Let $C=\{x\in U(A,\mathfrak{g})\mid x\cdot n\in N,\,\forall n\in N\}$. Then $C$ is an associative subalgebra of $U(A,\mathfrak{g})$. Since $N$ is an $A$-submodule, $\iota_A(A)\subseteq C$. Moreover, $\iota_{\mathfrak{g}}(x)\cdot n=\phi(x)(n)\in N$ for all $x\in\mathfrak{g}$ and $n\in N$, so $\iota_{\mathfrak{g}}(\mathfrak{g})\subseteq C$. Consequently, $C=U(A,\mathfrak{g})$. Hence $N$ is a $U(A,\mathfrak{g})$-submodule.
\end{proof}
\subsection{Filtrations and Graded Associative Algebras}
 We now introduce filtrations and graded associative algebras for Lie-Rinehart color algebras.

Let $(A,\mathfrak{g},\rho,\Gamma,(\cdot,\cdot))$ be a Lie-Rinehart color algebra and let $U(A,\mathfrak{g})$ be its universal enveloping algebra equipped with the natural $\Gamma$-grading structure.

For all $n\in\mathbb{Z}$, define the left $A$-submodule sequence of $U(A,\mathfrak{g})$:
If $n<0$, then $U_n=0$;
If $n=0$, then $U_0=A$;
If $n\ge 1$, then $U_n$ is the left $A$-submodule of $U(A,\mathfrak{g})$ generated by products of at most $n$ $\Gamma$-homogeneous elements in $\mathfrak{g}$.

The sequence $\{U_n\}_{n\in\mathbb{Z}}$ is an increasing multiplicative filtration of $U(A,\mathfrak{g})$, that is, for all $n,m\in\mathbb{Z}$,
\[
U_nU_m\subseteq U_{n+m}.
\]
Furthermore, each submodule $U_n$ is a $\Gamma$-graded submodule of $U(A,\mathfrak{g})$.

Define the quotient space $$\operatorname{gr}_n(U(A,\mathfrak{g}))=U_n(A,\mathfrak{g})/U_{n-1}(A,\mathfrak{g}).$$ The associated graded space is
\[
\operatorname{gr}(U(A,\mathfrak{g}))=\bigoplus_{n\ge 0}\operatorname{gr}_n(U(A,\mathfrak{g}))=\bigoplus_{n\ge 0}U_n(A,\mathfrak{g})/U_{n-1}(A,\mathfrak{g}).
\]
 $\operatorname{gr}(U(A,\mathfrak{g}))$ is a $\mathbb{Z}\times\Gamma$-graded $A$-algebra. The associative multiplication of $U(A,\mathfrak{g})$ naturally induces a graded multiplication on $\operatorname{gr}(U(A,\mathfrak{g}))$, such that 
\[
\operatorname{gr}_n(U(A,\mathfrak{g}))\cdot \operatorname{gr}_m(U(A,\mathfrak{g})) \subseteq \operatorname{gr}_{n+m}(U(A,\mathfrak{g})).
\]
for any non-negative integers $n,m$.
Take arbitrary $u\in (U_p)_\alpha$, $v\in (U_q)_\beta$, where $\alpha,\beta\in\Gamma$, $p,q\ge 0$. Then
\[
uv - (\alpha,\beta)vu \in U_{p+q-1}.
\]
It follows that in the associated graded algebra $\operatorname{gr}(U(A,\mathfrak{g}))$,
\[
\overline{u}\overline{v} = (\alpha,\beta)\overline{v}\overline{u},
\]
Then $\operatorname{gr}(U(A,\mathfrak{g}))$ is a color-commutative $\mathbb{Z}\times\Gamma$-graded $A$-algebra.

We use the color symmetric relations in \cite{ref22} to define
\(S_A(\mathfrak g)\). 
Let $J_A(\mathfrak{g})$ be the two-sided ideal of the tensor algebra $T_A(\mathfrak{g})$ generated by tensors of the form
\[
u\otimes v - (\alpha,\beta)v\otimes u,\quad u\in\mathfrak{g}_\alpha,\,v\in\mathfrak{g}_\beta,\,\alpha,\beta\in\Gamma.
\]
Define
\[
S_A(\mathfrak{g}):= T_A(\mathfrak{g})/J_A(\mathfrak{g}).
\]
Then $S_A(\mathfrak{g})$ is naturally equipped with a $\mathbb{Z}\times\Gamma$-grading. By the universal property of symmetric algebras, there exists a canonical $A$-color algebra homomorphism
\[
\varphi:S_A(\mathfrak{g})\to \operatorname{gr}(U(A,\mathfrak{g})).
\]
In next section, we will prove that $\varphi$ is an isomorphism of $A$-color algebras.
\subsection{Poincar\'e-Birkhoff-Witt Theroem}
\begin{lemma}
Let $(A,\mathfrak{g},\Gamma,(\cdot,\cdot),\rho)$ be a Lie--Rinehart color algebra and $B$ a color-commutative algebra.
Let $f:A\to B$ be a color algebra homomorphism, and $g:\mathfrak{g}\to B$ be a degree-preserving $\mathbb{R}$-linear map.
Suppose that for all homogeneous elements $a\in A$ and $x\in\mathfrak{g}$, we have $g(ax)=f(a)g(x)$.
Then there exists a unique graded associative algebra homomorphism $\varphi:S_A(\mathfrak{g})\to B$ such that $\varphi\circ j_A = f$ and $\varphi\circ j_\mathfrak{g}=g$, where $j_A:A\to S_A(\mathfrak{g})$, $j_\mathfrak{g}:\mathfrak{g}\to S_A(\mathfrak{g})$ are the canonical maps.
\end{lemma}

\begin{proof}
For homogeneous \(a\in A_\alpha\) and \(x\in\mathfrak g_\beta\),
we have
\[
g(x)f(a)=(\beta,\alpha)f(a)g(x).
\]
Hence, by the universal property of \(T_A(\mathfrak g)\),
there exists a unique graded associative algebra homomorphism
\(\widetilde{\varphi}:T_A(\mathfrak g)\to B\) satisfying
\[
\widetilde{\varphi}\circ\widetilde{j}_A=f,
\qquad
\widetilde{\varphi}\circ\widetilde{j}_{\mathfrak{g}}=g,
\]
where $\widetilde{j_A}:A\to T_A(\mathfrak{g})$ and $\widetilde{j_\mathfrak{g}}:\mathfrak{g}\to T_A(\mathfrak{g})$ denote the canonical maps. 
Take arbitrary homogeneous elements $x\in\mathfrak{g}_\alpha$, $y\in\mathfrak{g}_\beta$.
\[
\widetilde{\varphi}(x\otimes_A y - (\alpha,\beta)y\otimes_A x)
= g(x)g(y)-(\alpha,\beta)g(y)g(x)=0.
\]
Hence $J_A(\mathfrak{g})\subseteq \ker\widetilde{\varphi}$.
Let $\pi:T_A(\mathfrak{g})\to S_A(\mathfrak{g})=T_A(\mathfrak{g})/J_A(\mathfrak{g})$.
By properties of quotient algebras, there exists a unique color algebra homomorphism $\varphi:S_A(\mathfrak{g})\to B$ such that $\varphi\circ\pi=\widetilde{\varphi}$.
We have
\[
\varphi\circ j_A = \varphi\circ\pi\circ \widetilde{j_A} = \widetilde{\varphi}\circ \widetilde{j_A} = f,
\]
\[
\varphi\circ j_\mathfrak{g} = \varphi\circ\pi\circ \widetilde{j_\mathfrak{g}} = \widetilde{\varphi}\circ \widetilde{j_\mathfrak{g}}= g.
\]
\end{proof}
\begin{theorem}\label{1349}
Let $(A,\mathfrak{g},\rho,(\cdot,\cdot),\Gamma)$ be a Lie-Rinehart color algebra. Assume that $A$ is a free $\mathbb{R}$-module and  $\mathfrak g$ is a free
$A$-module, with homogeneous $\mathbb{R}$ and $A$ bases, respectively. Then $\varphi:S_A(\mathfrak{g})\to \operatorname{gr}(U(A,\mathfrak{g}))$ is an isomorphism of $A$-color algebras.
\end{theorem}

\begin{proof}
Define the map $\bar{\iota}_A:A\to \operatorname{gr}U(A,\mathfrak{g})$ by $\bar{\iota}_A(a)=\iota_A(a)+U_{-1}$. Then $\bar{\iota}_A$ is color algebra homomorphism. Define another map $\bar{\iota}_\mathfrak{g}:\mathfrak{g}\to \operatorname{gr}U(A,\mathfrak{g})$ by $\bar{\iota}_\mathfrak{g}(x)=\iota_\mathfrak{g}(x)+U_0$. We have
\[
\bar{\iota}_\mathfrak{g}(ax)=\iota_\mathfrak{g}(ax)+U_0=\iota_A(a)\iota_\mathfrak{g}(x)+U_0=\big(\iota_A(a)+U_{-1}\big)\big(\iota_\mathfrak{g}(x)+U_0\big)=\bar{\iota}_A(a)\bar{\iota}_\mathfrak{g}(x).
\]
\[
\bar{\iota}_\mathfrak{g}(x)\bar{\iota}_\mathfrak{g}(y)=(\beta,\gamma)\bar{\iota}_\mathfrak{g}(y)\bar{\iota}_\mathfrak{g}(x).
\]
 for homogeneous elements $a\in A_\alpha$, $x\in\mathfrak{g}_\beta$, $y\in\mathfrak{g}_\gamma$.
By the universal property of $S_A(\mathfrak{g})$, there exists a unique graded $A$-color algebra homomorphism $$\varphi:S_A(\mathfrak{g})\to \operatorname{gr}U(A,\mathfrak{g})$$ such that $\varphi(a)=\bar{\iota}_A(a)$, $\varphi(x)=\bar{\iota}_\mathfrak{g}(x)=\iota_\mathfrak{g}(x)+U_0$ for all $x\in\mathfrak{g}$, $a\in A$. Next We prove that $\varphi$ is an isomorphism.

Let $\mathfrak g_1$ be the underlying graded $A$-module of
$\mathfrak g$ equipped with the zero bracket and the zero anchor
$\rho_1$. Then $(A,\mathfrak g_1)$ is a Lie-Rinehart color algebra.

Set $W=U^+(A\rtimes\mathfrak{g})$. Let $J$ be the two-sided ideal of $W$ generated by the set $$\{a\otimes z - a\cdot z\mid a\in A,\,z\in A\rtimes\mathfrak{g}\},$$ where $a\cdot z$ denotes the left $A$-module action on $A\rtimes\mathfrak{g}$. Then $U(A,\mathfrak{g})\cong W/J$. Similarly, set $W'=U^+(A\rtimes\mathfrak{g}_1)$, and let $I$ be the two-sided ideal of $W'$ generated by $$\{a\otimes z - a\cdot z\mid a\in A,\,z\in A\rtimes\mathfrak{g}_1\}$$ where $a\cdot z$ denotes the left $A$-module action on $A\rtimes\mathfrak{g}_1$. Then $U(A,\mathfrak{g}_1)\cong W'/I$.

Let $\{W_p\}_{p\ge0}$ be the natural filtration on $W$ determined by the number of $\mathfrak{g}$-factors, and let $\{W_p'\}_{p\ge0}$ be the natural filtration on $W'$ determined by the number of $\mathfrak{g}_1$-factors. Define $U_p=(W_p+J)/J$. By the isomorphism theorem for quotient spaces,
\[
U_p/U_{p-1}\cong \frac{W_p+J}{W_{p-1}+J}\cong \frac{W_p}{W_p\cap(W_{p-1}+J)}=\frac{W_p}{W_{p-1}+(W_p\cap J)}.
\]
Set $H_p=W_{p-1}+(W_p\cap J)$. Then $$\operatorname{gr}U(A,\mathfrak{g})\cong \bigoplus_{p\ge0}W_p/H_p.$$ Similarly, set $H_p'=W_{p-1}'+(W_p'\cap I)$. Then $$\operatorname{gr}U(A,\mathfrak{g}_1)\cong \bigoplus_{p\ge0}W_p'/H_p'.$$

Let $\{z_\lambda\mid\lambda\in\Lambda_A\}$ be a homogeneous $\mathbb{R}$-basis of $A$, and $\{z_\lambda\mid\lambda\in\Lambda_\mathfrak{g}\}$ a homogeneous $\mathbb{R}$-basis of $\mathfrak{g}$, with $\Lambda_A\cap\Lambda_\mathfrak{g}=\emptyset$. Let $\Lambda=\Lambda_A\cup\Lambda_\mathfrak{g}$, fix a total order on $\Lambda$ such that $\lambda<\mu$ whenever $\lambda\in\Lambda_A$ and $\mu\in\Lambda_\mathfrak{g}$. Denote by $\{z_\lambda'\mid\lambda\in\Lambda\}$ the corresponding homogeneous basis of $A\oplus\mathfrak{g}_1$.

By the PBW theorem for Lie color algebras, $(z_{\lambda_1}')^{r_1}\cdots(z_{\lambda_n}')^{r_n}$ forms an $\mathbb{R}$-basis of $W'$, and $z_{\lambda_1}^{r_1}\cdots z_{\lambda_n}^{r_n}$ forms an $\mathbb{R}$-basis of $W$, where $n\ge0$, $\lambda_1<\cdots<\lambda_n$, $r_i\in\mathbb{N}$. If $(|z_i|,|z_i|)=1$, then $r_i\in\mathbb{N}$; if $(|z_i|,|z_i|)=-1$, then $r_i\in\{0,1\}$.

Define an $\mathbb{R}$-linear bijection $\Psi:W'\to W$ such that $$\Psi\big((z_{\lambda_1}')^{r_1}\cdots(z_{\lambda_n}')^{r_n}\big)=z_{\lambda_1}^{r_1}\cdots z_{\lambda_n}^{r_n}.$$ Clearly $\Psi$ is $\mathbb{R}$-linear and bijective. In particular, $\Psi$ preserves filtrations and $\Gamma$-gradings, i.e., $\Psi\big((W_p')_\alpha\big)=(W_p)_\alpha$ for all $p\ge0$ and $\alpha\in\Gamma$.

Take arbitrary homogeneous monomials $u\in(W_p')_\alpha$, $v\in(W_q')_\beta$. From the color commutation relations in $W,W'$, we obtain
\[
\Psi(uv)-\Psi(u)\Psi(v)\in(W_{p+q-1})_{\alpha+\beta}.
\]

We next show that $\Psi(H_p')=H_p$. For homogeneous elements $a\in A$ and $z\in A\rtimes\mathfrak{g}_1$, write $r'(a,z)=a\otimes z - a\cdot z\in W'$. Let $r(a,z)\in W$ be the element corresponding to $r'(a,z)$.

Take any $h=h_0+h_1\in H_p'$, $h_0\in W_{p-1}'$, $h_1\in W_p'\cap I$. Since $I$ is generated by all $r'(a,z)$, $W'$ is graded by the number of $\mathfrak g_1$-factors and
the generators of $I$ are homogeneous for this grading, then we may
retain only the summands of degree $p$ modulo $W'_{p-1}$. Here
\[
h_1\equiv \sum_{k=1}^n u_k'r'(a_k,z_k)v_k'\pmod{W_{p-1}'},
\]
where $u_k',v_k'\in W_p'$, and the total filtration degree of each term $u_k'r'(a_k,z_k)v_k'$ is at most $p$.

By the product formula above, we have
\[
\Psi\big(u_k'r'(a_k,z_k)v_k'\big)-\Psi(u_k')r(a_k,z_k)\Psi(v_k')\in W_{p-1}.
\]
Thus $\Psi(h_1)\in J+W_{p-1}$. Since $\Psi(h_1)\in W_p$, we get
\[
\Psi(h_1)\in (W_p\cap J)+W_{p-1}=H_p.
\]
Since $\Psi(h_0)\in W_{p-1}\subseteq H_p$, we have $\Psi(h)=\Psi(h_0)+\Psi(h_1)\in H_p$. Hence $\Psi(H_p')\subseteq H_p$.

To prove the reverse inclusion, fix a totally ordered homogeneous
\(A\)-basis \(\{x_i\}_{i\in I}\) of \(\mathfrak g\) and the
corresponding basis of \(\mathfrak g_1\). In \(W/J\), we have
\[
\begin{aligned}
x_i a
  &= (|x_i|,|a|)\,a x_i+\rho(x_i)(a),\\
x_i x_j
  &= (|x_i|,|x_j|)\,x_jx_i+[x_i,x_j],
     && i>j,\\
x_i^2
  &=\frac12[x_i,x_i],
     && (|x_i|,|x_i|)=-1,
\end{aligned}
\]
where
\[
[x_i,x_j]=\sum_{k\in I}c_{ij}^{\,k}x_k,
\qquad c_{ij}^{\,k}\in A.
\]
Consequently, every element of \(W/J\) can be written as a finite
\(A\)-linear combination of ordered monomials
\[
\sum_{i_1\leq\cdots\leq i_r}
a_{i_1,\ldots,i_r}x_{i_1}\cdots x_{i_r},
\qquad a_{i_1,\ldots,i_r}\in A,
\]
where an index \(i\) with \((|x_i|,|x_i|)=-1\) occurs at most once
in each monomial. We show that the resulting expression is
independent of the order in which the above relations are applied.
For a monomial \(M\), let \(n(M)\), \(i(M)\), and \(c(M)\)
be the numbers of basis factors, inversions, and coefficient
factors, respectively. For each monomial \(M'\) obtained from
\(M\) by applying a relation, we have
\[
\begin{aligned}
&n(M')<n(M),\quad\text{or}\\
&n(M')=n(M),\quad i(M')<i(M),\quad\text{or}\\
&n(M')=n(M),\quad i(M')=i(M),\quad c(M')<c(M).
\end{aligned}
\]
Thus the process ends after finitely many steps.

For \(x_i x_j a\), we may first interchange \(x_i\) and \(x_j\),
or first move \(a\) to the left of \(x_j\). Direct calculations
for the following monomials show that both orders of rearrangement
give the same result:
\[
abc,\quad x_iab,\quad x_ix_ja,\quad x_i^2a,\quad
x_ix_jx_k,\quad x_i^2x_j,\quad x_ix_j^2,\quad x_i^3.
\]

Thus induction gives a unique ordered expression.

For $u,v\in W$, the terms $u(a\otimes z)v$ and $u(a\cdot z)v$
have the same ordered expression. Thus every element of $J$ has
zero ordered expression. The same argument with zero bracket and
anchor gives unique ordered expressions in $W'$ and shows that
every element of $I$ has zero ordered expression.

Since $\Psi^{-1}(W_{p-1})\subseteq W'_{p-1}$, it suffices to show
that $\Psi^{-1}(j)\in H'_p$ for every $j\in W_p\cap J$.
The degree-$p$ parts of the ordered expressions of $j$ and
$\Psi^{-1}(j)$ agree under the correspondence of basis factors,
because the bracket and anchor terms lower the degree.
Since $j\in J$, both degree-$p$ parts vanish.

Let $j'_p$ be the degree-$p$ component of $\Psi^{-1}(j)$ in $W'$.
Rearrangement in $W'$ preserves degree and gives equalities modulo $I$,
so the vanishing of the ordered expression of $j'_p$ gives
$j'_p\in W'_p\cap I$. Therefore
\[
\Psi^{-1}(j)\in j'_p+W'_{p-1}\subseteq H'_p,
\]
which proves $\Psi^{-1}(H_p)\subseteq H'_p$.
Hence,
\[
\Psi(H'_p)=H_p.
\]

Define the linear map $\bar{\Psi}_p:W_p'/H_p'\to W_p/H_p$ by $\bar{\Psi}_p(u+H_p')=\Psi(u)+H_p$ for all $u\in W_p'$. Since $\Psi$ is bijective and $\Psi(H_p')=H_p$, $\bar{\Psi}_p$ is a linear isomorphism. Let $\bar{\Psi}=\bigoplus_{p\ge0}\bar{\Psi}_p$,
\[
\bar{\Psi}:\operatorname{gr}U(A,\mathfrak{g}_1)\to \operatorname{gr}U(A,\mathfrak{g}).
\]
Let $u\in (W'_p)_\alpha$ and $v\in (W'_q)_\beta$. Since
\[
\Psi(uv)-\Psi(u)\Psi(v)
\in (W_{p+q-1})_{\alpha+\beta}
\subseteq (H_{p+q})_{\alpha+\beta},
\]
we have
\[
\Psi(uv)+H_{p+q}=\Psi(u)\Psi(v)+H_{p+q}.
\]
It follows that
\[
\begin{aligned}
\overline{\Psi}\bigl((u+H'_p)(v+H'_q)\bigr)
&=\overline{\Psi}(uv+H'_{p+q})\\
&=\Psi(u)\Psi(v)+H_{p+q}\\
&=\overline{\Psi}(u+H'_p)\,\overline{\Psi}(v+H'_q).
\end{aligned}
\]
Moreover, $\Psi$ preserves the $\Gamma$-grading, so
$\overline{\Psi}$ preserves the $\mathbb Z_{\ge 0}\times\Gamma$-bigrading.
Hence $\overline{\Psi}$ is an isomorphism of algebras.
Since $\mathfrak g_1$ has zero bracket and zero anchor,
$U(A,\mathfrak g_1)$ is color-commutative. The canonical maps satisfy
$\iota_{\mathfrak g_1}(ax)=\iota_{A_1}(a)\iota_{\mathfrak g_1}(x)$.
Identifying $\mathfrak g_1$ with $\mathfrak g$ as graded $A$-modules,
we obtain by Lemma~4.12 a unique homomorphism of $A$-color algebras
\[
\theta_0\colon S_A(\mathfrak g)\longrightarrow U(A,\mathfrak g_1)
\]
such that $\theta_0\circ j_A=\iota_{A_1}$ and
$\theta_0\circ j_{\mathfrak g}=\iota_{\mathfrak g_1}$,
where $j_A$ and $j_{\mathfrak g}$ are the canonical maps into
$S_A(\mathfrak g)$.

Conversely, the color commutativity of $S_A(\mathfrak g)$ and the
$A$-linearity of $j_{\mathfrak g}$ ensure that $j_A$ and $j_{\mathfrak g}$
satisfy the compatibility conditions for $(A,\mathfrak g_1)$.
By the universal property of $U(A,\mathfrak g_1)$, there is a unique
graded algebra homomorphism
$\eta\colon U(A,\mathfrak g_1)\to S_A(\mathfrak g)$ such that
$\eta\circ\iota_{A_1}=j_A$ and
$\eta\circ\iota_{\mathfrak g_1}=j_{\mathfrak g}$.

Since
\[
(\eta\circ\theta_0)\circ j_A=j_A,\qquad (\eta\circ\theta_0)\circ j_\mathfrak{g}=j_\mathfrak{g},
\]
\[
(\theta_0\circ\eta)\circ\iota_{A_1}=\iota_{A_1},\qquad (\theta_0\circ\eta)\circ\iota_{\mathfrak{g}_1}=\iota_{\mathfrak{g}_1},
\]
we have $\eta\circ\theta_0=\operatorname{id}_{S_A(\mathfrak{g})}$, $\theta_0\circ\eta=\operatorname{id}_{U(A,\mathfrak{g}_1)}$, so $\theta_0$ is an isomorphism.

Since $U(A,\mathfrak{g}_1)$ is a $\mathbb{Z}_{\ge0}$-graded algebra, there exists a natural isomorphism $\kappa:U(A,\mathfrak{g}_1)\to \operatorname{gr}U(A,\mathfrak{g}_1)$. Therefore, $\theta=\kappa\circ\theta_0:S_A(\mathfrak{g})\to\operatorname{gr}U(A,\mathfrak{g}_1)$ is an isomorphism. Since $\varphi=\bar{\Psi}\circ\theta$, $\varphi$ is an isomorphism of $A$-color algebras.
\end{proof}

\begin{corollary}\label{1348}
The map $\iota_\mathfrak{g}:\mathfrak{g}\to U(A,\mathfrak{g})$ is injective.
\end{corollary}
\begin{proof}
From the above theorem, the map $$\mathfrak{g}= S^1_A(\mathfrak{g})\to U_1/U_0, x\mapsto \iota_\mathfrak{g}(x)+U_0$$ is bijective. If $\iota_\mathfrak{g}(x)=0$, then $\iota_\mathfrak{g}(x)+U_0=0$. Hence $\ker\iota_\mathfrak{g}=0$, so $\iota_\mathfrak{g}$ is injective.
\end{proof}
The following Poincar\'e--Birkhoff--Witt theorem for Lie-Rinehart color algebras follows directly from Theorem \ref{1349} and Corollory \ref{1348}.
\begin{theorem}
Let $\{x_\lambda\mid\lambda\in\Lambda\}$ be a homogeneous $A$-basis of $\mathfrak{g}$, where $\Lambda$ is a well-ordered set and $x_\lambda\in\mathfrak{g}_{\gamma_\lambda}$. Then $U(A,\mathfrak g)$ is generated as a left $A$-module by
the monomials
\[
\iota_{\mathfrak g}(x_{\lambda_1})\cdots
\iota_{\mathfrak g}(x_{\lambda_n}),
\qquad n\geq 0,\qquad
\lambda_1\leq\cdots\leq\lambda_n,
\]
where $\lambda_j<\lambda_{j+1}$ whenever $1\leq j<n$ and
$(\gamma_{\lambda_j},\gamma_{\lambda_j})=-1$.
For $n=0$, the monomial is the unit of $U(A,\mathfrak g)$.

\end{theorem}
\subsection*{Acknowledgements}
The authors are supported by the  NSF of Hei Longjiang Province (No. LH2024A014).
\subsection*{Declaration of competing interest}
There is no competing interest.


\end{document}